\documentclass[journal]{IEEEtran}
\usepackage[latin9]{inputenc}
\usepackage{verbatim}
\usepackage{amsthm}
\usepackage{amsmath}
\usepackage{amssymb}
\usepackage{graphicx}
\usepackage{graphics}
\usepackage{color}

\makeatletter
\theoremstyle{plain}
\newtheorem{thm}{\protect\theoremname}
\theoremstyle{remark}
\newtheorem{rem}[thm]{\protect\remarkname}
\theoremstyle{definition}
\newtheorem{problem}[thm]{\protect\problemname}
\theoremstyle{definition}
\newtheorem{defn}[thm]{\protect\definitionname}
\theoremstyle{plain}
\newtheorem{lem}[thm]{\protect\lemmaname}
\theoremstyle{plain}
\newtheorem{cor}[thm]{\protect\corollaryname}
\theoremstyle{plain}
\newtheorem{prop}[thm]{\protect\propositionname}

\usepackage{amsfonts}
\usepackage{amsmath}
\usepackage{amssymb}
\usepackage{epsfig}
\usepackage{graphics}
\usepackage{graphicx}
\usepackage{needspace}
\usepackage{color, colortbl}
\definecolor{Gray}{gray}{0.9}
\epsfclipon

\newtheorem{proposition}{Proposition}

\makeatother

\providecommand{\lemmaname}{Lemma}
\providecommand{\problemname}{Problem}
\providecommand{\propositionname}{Proposition}
\providecommand{\theoremname}{Theorem}

\makeatother

\providecommand{\corollaryname}{Corollary}
\providecommand{\definitionname}{Definition}
\providecommand{\lemmaname}{Lemma}
\providecommand{\problemname}{Problem}
\providecommand{\propositionname}{Proposition}
\providecommand{\remarkname}{Remark}
\providecommand{\theoremname}{Theorem}

\begin{document}

\title{Optimal Resource Allocation for Network Protection Against Spreading Processes}

\author{Victor M. Preciado, Michael Zargham, Chinwendu Enyioha, Ali Jadbabaie,
and George Pappas %
\thanks{The authors are with the Department of Electrical and Systems Engineering
at the University of Pennsylvania, Philadelphia PA 19104.
}}
\maketitle
\begin{abstract}
We study the problem of containing spreading processes in arbitrary
directed networks by distributing protection resources throughout
the nodes of the network. We consider two types of protection resources
are available: (\emph{i}) \emph{Preventive} resources able to defend
nodes against the spreading (such as vaccines in a viral infection
process), and (\emph{ii}) \emph{corrective} resources able to neutralize
the spreading after it has reached a node (such as antidotes). We
assume that both preventive and corrective resources have an associated
cost and study the problem of finding the cost-optimal distribution
of resources throughout the nodes of the network. We analyze these
questions in the context of viral spreading processes in directed
networks. We study the following two problems: (\emph{i}) Given a
fixed budget, find the optimal allocation of preventive and corrective
resources in the network to achieve the highest level of containment,
and (\emph{ii}) when a budget is not specified, find the minimum budget
required to control the spreading process. We show that both resource
allocation problems can be solved in polynomial time using Geometric
Programming (GP) for arbitrary \emph{directed }graphs of \emph{nonidentical}
nodes and a wide class of cost functions. We illustrate our approach by designing optimal protection strategies
to contain an epidemic outbreak that propagates through an air transportation
network.
\end{abstract}

\section{Introduction}

Understanding spreading processes in complex networks and designing
control strategies to contain them are relevant problems in many different
settings, such as epidemiology and public health \cite{Bai75}, computer
viruses \cite{GGT03}, or security of cyberphysical networks \cite{roy2012security}.
In this paper, we analyze the problem of controlling spreading processes
in networks by distributing protection resources throughout the nodes.
In our study, we consider two types of containment resources: (\emph{i})
\emph{Preventive} resources able to protect (or `immunize') nodes
against the spreading (such as vaccines in a viral infection process),
and (\emph{ii}) \emph{corrective} resources able to neutralize the
spreading after it has reached a node, such as antidotes in a viral
infection. In our framework, we associate a cost with these resources
and study the problem of finding the cost-optimal distribution of
resources throughout the network to contain the spreading.

In the literature, we find several approaches to model spreading mechanisms
in arbitrary contact networks. The analysis
of this question in arbitrary (undirected) contact networks was first studied by
Wang et al. \cite{WCWF03} for a Susceptible-Infected-Susceptible
(SIS) discrete-time model. In \cite{GMT05}, Ganesh et al. studied
the epidemic threshold in a continuous-time SIS spreading processes.
In both continuous- and discrete-time models, there is a close connection
between the speed of the spreading and the spectral radius of the
network (i.e., the largest eigenvalue of its adjacency matrix) \cite{MOK09}. In contrast to most current research, we focus our attention on the analysis of directed contact networks. From a practical point of view, many real networks are more naturally modeled using weighted and directed edges.

Designing strategies to contain spreading processes in networks is
a central problem in public health and network security. In this context,
the following question is of particular interest: given a contact
network (possibly weighted and/or directed) and resources that provide
partial protection (e.g., vaccines and/or antidotes), how should one
distribute these resources throughout the networks in a cost-optimal
manner to contain the spread? This question has been addressed
in several papers. Cohen et al. \cite{cohen2003efficient} proposed
a heuristic vaccination strategy called \emph{acquaintance immunization
policy} and proved it to be much more efficient than random vaccine
allocation. In \cite{BCGS10}, Borgs et al. studied theoretical
limits in the control of spreads in undirected network with a non-homogeneous
distribution of antidotes. Chung et al. \cite{chung2009distributing}
studied a heuristic immunization strategy based on the PageRank vector
of the contact graph. In the control systems literature, Wan et al.
proposed in \cite{WRS08} a method to design control strategies in
undirected networks using eigenvalue sensitivity analysis. Our work
is related to that in \cite{GOM11}, where the authors study the problem
of minimizing the level of infection in an undirected network using
corrective resources within a given budget. In \cite{PDS13} a linear-fractional optimization
program was proposed to compute the optimal investment on disease
awareness over the nodes of a social network to contain a spreading
process. Also, in \cite{PZ13,PZEJP13},
the authors proposed a convex formulation to find the optimal allocation
of protective resources in an undirected network using semidefinite
programming (SDP).

Breaking the network symmetry prevents us from using previously proposed approaches to find the optimal resource allocation for network protection. In this paper, we propose a novel formulation based on geometric programming to find the optimal allocation of protection resources in \emph{weighted and directed} networks
of \emph{nonidentical} agents in polynomial time.

The paper is organized as follows. In Section \ref{Notation}, we
introduce notation and background needed in our derivations. We state
the resource allocation problems solved in this paper in Subsection
\ref{sub:Problem-Statements}. In Section \ref{sec:Convex Framework},
we propose a convex optimization framework to efficiently solve the
allocation problems in polynomial time. Subsection \ref{sub:GP for SC Digraphs},
present the solution to the allocation problem for strongly connected
graphs. We extend this result to general directed graphs (not necessarily
strongly connected) in Subsection \ref{sub:Allocation-for general digraphs}.
We illustrate our results using a real-world air transportation network
in Section \ref{sub:Numerical-Results}. We include some conclusions
in Section \ref{sec:Conclusions}.

\section{\label{Notation}Preliminaries \& Problem Definition}

We introduce notation and preliminary results needed in our derivations.
In the rest of the paper, we denote by $\mathbb{R}_{+}^{n}$ (respectively,
$\mathbb{R}_{++}^{n}$) the set of $n$-dimensional vectors with nonnegative
(respectively, positive) entries. We denote vectors using boldface letters
and matrices using capital letters. $I$ denotes the identity matrix
and $\mathbf{1}$ the vector of all ones. $\Re\left(z\right)$ denotes
the real part of $z\in\mathbb{C}$.

\subsection{\label{sub:Graph-Theory}Graph Theory}

A \emph{weighted}, \emph{directed} graph (also called digraph) is
defined as the triad $\mathcal{G}\triangleq\left(\mathcal{V},\mathcal{E},\mathcal{W}\right)$,
where (\emph{i}) $\mathcal{V}\triangleq\left\{ v_{1},\dots,v_{n}\right\} $
is a set of $n$ nodes, (\emph{ii}) $\mathcal{E}\subseteq\mathcal{V}\times\mathcal{V}$
is a set of ordered pairs of nodes called directed edges, and (\emph{iii})
the function $\mathcal{W}:\mathcal{E}\rightarrow\mathbb{R}_{++}$
associates \textit{positive} real weights to the edges in $\mathcal{E}$.
By convention, we say that $\left(v_{j},v_{i}\right)$ is an edge
from $v_{j}$ pointing towards $v_{i}$. We define the in-neighborhood
of node $v_{i}$ as $\mathcal{N}_{i}^{in}\triangleq\left\{ j:\left(v_{j},v_{i}\right)\in\mathcal{E}\right\} $,
i.e., the set of nodes with edges pointing towards $v_{i}$. We define
the weighted \emph{in-degree} (resp., out-degree) of node $v_{i}$
as $\deg_{in}\left(v_{i}\right)\triangleq\sum_{j\in\mathcal{N}_{i}^{in}}\mathcal{W}\left(\left(v_{j},v_{i}\right)\right)$
(resp., $\deg_{out}\left(v_{i}\right)\triangleq\sum_{j\in\mathcal{N}_{i}^{out}}\mathcal{W}\left(\left(v_{j},v_{i}\right)\right)$).
A directed path from $v_{i_{1}}$ to $v_{i_{l}}$ in $\mathcal{G}$
is an ordered set of vertices $\left(v_{i_{1}},v_{i_{2}},\ldots,v_{i_{l+1}}\right)$
such that $\left(v_{i_{s}},v_{i_{s+1}}\right)\in\mathcal{E}$ for
$s=1,\ldots,l$. A directed graph $\mathcal{G}$ is \emph{strongly
connected} if, for every pair of nodes $v_{i},v_{j}\in\mathcal{V}$,
there is a directed path from $v_{i}$ to $v_{j}$.

The \emph{adjacency matrix} of a weighted, directed graph $\mathcal{G}$,
denoted by $A_{\mathcal{G}}=[a_{ij}]$, is an $n\times n$ matrix defined
entry-wise as $a_{ij}=\mathcal{W}((v_{j},v_{i}))$ if edge $(v_{j},v_{i})\in\mathcal{E}$,
and $a_{ij}=0$ otherwise. Given an $n\times n$ matrix $M$, we denote
by $\mathbf{v}_{1}\left(M\right),\ldots,\mathbf{v}_{n}\left(M\right)$
and $\lambda_{1}\left(M\right),\ldots,\lambda_{n}\left(M\right)$
the set of eigenvectors and corresponding eigenvalues of $M$, respectively,
where we order them in decreasing order of their real parts, i.e.,
$\Re\left(\lambda_{1}\right)\geq\Re\left(\lambda_{2}\right)\geq\ldots\geq\Re\left(\lambda_{n}\right)$.
We respectively call $\lambda_{1}\left(M\right)$ and $\mathbf{v}_{1}\left(M\right)$
the dominant eigenvalue and eigenvector of $M$. The spectral radius
of $M$, denoted by $\rho\left(M\right)$, is the maximum modulus
across all eigenvalue of $M$.

In this paper, we only consider graphs with positively weighted edges;
hence, the adjacency matrix of a graph is always nonnegative. Conversely,
given a $n\times n$ nonnegative matrix $A$, we can associate a directed
graph $\mathcal{G}_{A}$ such that $A$ is the adjacency matrix of
$\mathcal{G}_{A}$. Finally, a nonnegative matrix $A$ is \emph{irreducible}
if and only if its associated graph $\mathcal{G}_{A}$ is strongly
connected.

\textbf{}%

\subsection{\label{sub:Epidemic-Model}Stochastic Spreading Model in Arbitrary
Networks}

A popular stochastic model to simulate spreading processes is the
so-called susceptible-infected-susceptible (SIS) epidemic model, first
introduced by Weiss and Dishon \cite{weiss1971asymptotic}. Wang et
al. \cite{WCWF03} proposed a discrete-time extension of the SIS model
to simulate spreading processes in networked populations. A continuous-time
version, called the N-intertwined SIS model, was recently proposed
and rigorously analyzed by Van Mieghem et al. in \cite{MOK09}. In
this paper, we formulate our problem using a further extension of
the SIS model recently proposed in \cite{VO13}. We call this model
the \emph{Neterogeneous Networked SIS model} (HeNeSIS).

This HeNeSIS model is a continuous-time networked Markov process in
which each node in the network can be in one out of two possible states,
namely, susceptible or infected. Over time, each node $v_{i}\in\mathcal{V}$
can change its state according to a stochastic process parameterized
by (\emph{i}) the node infection rate $\beta_{i}$, and (\emph{ii})
its recovery rate $\delta_{i}$. In our work, we assume that both
$\beta_{i}$ and $\delta_{i}$ are node-dependent and adjustable via
the injection of vaccines and/or antidotes in node $v_{i}$.

The evolution of the HeNeSIS model can be described as follows. The
state of node $v_{i}$ at time $t\geq0$ is a binary random variable
$X_{i}\left(t\right)\in\{0,1\}$. The state $X_{i}\left(t\right)=0$
(resp., $X_{i}\left(t\right)=1$) indicates that node $v_{i}$ is
in the susceptible (resp., infected) state. We define the vector of
states as $X\left(t\right)=\left(X_{1}\left(t\right),\ldots,X_{n}\left(t\right)\right)^{T}$.
The state of a node can experience two possible stochastic transitions:
\begin{enumerate}
\item Assume node $v_{i}$ is in the susceptible state at time $t$. This
node can switch to the infected state during the (small) time interval
$\left[t,t+\Delta t\right)$ with a probability that depends on: (\emph{i})
its infection rate $\beta_{i}>0$, (\emph{ii}) the strength of its
incoming connections $\left\{ a_{ij},\mbox{ for }j\in\mathcal{N}_{i}^{in}\right\} $,
and (\emph{iii}) the states of its in-neighbors $\left\{ X_{j}\left(t\right),\mbox{ for }j\in\mathcal{N}_{i}^{in}\right\} $.
Formally, the probability of this transition is given by 
\begin{multline}
\Pr\left(X_{i}(t+\Delta t)=1|X_{i}(t)=0,X(t)\right)=\\
\sum_{j\in\mathcal{N}_{i}^{in}}a_{ij}\beta_{i}X_{j}\left(t\right)\Delta t+o(\Delta t),
\end{multline}
where $\Delta t>0$ is considered an asymptotically small time interval.
\item Assuming node $v_{i}$ is infected, the probability of $v_{i}$ recovering
back to the susceptible state in the time interval $\left[t,t+\Delta t\right)$
is given by
\begin{equation}
\Pr(X_{i}(t+\Delta t)=0|X_{i}(t)=1,X(t))=\delta_{i}\Delta t+o(\Delta t),
\end{equation}
where $\delta_{i}>0$ is the curing rate of node $v_{i}$.
\end{enumerate}
This HeNeSIS model is therefore a continuous-time Markov process with
$2^{n}$ states in the limit $\Delta t\to0^{+}$. Unfortunately, the
exponentially increasing state space makes this model hard to analyze
for large-scale networks. To overcome this
limitation, we use a mean-field approximation of its dynamics \cite{barrat2008dynamical}.
This approximation is widely used in the field of epidemic analysis and control \cite{WCWF03}-\cite{MOK09},\cite{WRS08}-\cite{GMWPM12}, since it performs numerically well for many realistic network topologies\footnote{Finding rigorous conditions on the network structure for the mean-field approximation to be tight is a matter of current research in the community and beyond the scope of this paper. For a recent study on the accuracy of the mean-field approximation in realistic network topologies, see [19].}. Using the Kolmogorov forward equations and a mean-field approach,
one can approximate the dynamics of the spreading process using a
system of $n$ ordinary differential equations, as follows. Let us
define $p_{i}\left(t\right)\triangleq\Pr\left(X_{i}\left(t\right)=1\right)=E\left(X_{i}\left(t\right)\right)$,
i.e., the marginal probability of node $v_{i}$ being infected at
time $t$. Hence, the Markov differential equation \cite{van2006performance}
for the state $X_{i}\left(t\right)=1$ is the following,
\begin{equation}
\frac{dp_{i}\left(t\right)}{dt}=\left(1-p_{i}\left(t\right)\right)\beta_{i}\sum_{j=1}^{n}a_{ij}p_{j}\left(t\right)-\delta_{i}p_{i}\left(t\right).\label{eq:HeNiSIS dynamics}
\end{equation}
Considering $i=1,\ldots,n$, we obtain a system of nonlinear differential
equation with a complex dynamics. In the following, we derive a sufficient
condition for infections to die out exponentially fast.

We can write the mean-field approximation of the HeNeSIS model in matrix form as 
\begin{equation}
\frac{d\boldsymbol{p}\left(t\right)}{dt}=\left(BA_{\mathcal{G}}-D\right)\boldsymbol{p}\left(t\right)-P\left(t\right)BA_{\mathcal{G}}\boldsymbol{p}\left(t\right),\label{eq:heteroSIS}
\end{equation}
where $\boldsymbol{p}\left(t\right)\triangleq\left(p_{1}\left(t\right),\ldots,p_{n}\left(t\right)\right)^{T}$,
$B\triangleq\mbox{diag}(\beta_{i})$, $D\triangleq\mbox{diag}\left(\delta_{i}\right)$,
and $P\left(t\right)\triangleq\mbox{diag}(p_{i}\left(t\right))$.
This ODE presents an equilibrium point at $\boldsymbol{p}^{*}=0$,
called the disease-free equilibrium. In practice, the levels of infection in the population are very small, allowing us to linearize the dynamics around the disease-free equilibrium. Furthermore, it was proved in \cite{PZEJP13} that the linearized dynamics upper-bounds the nonlinear one. Therefore, we can stabilize the nonlinear dynamics by stabilizing its linear approximation. More specifically, a stability analysis of this
ODE around the equilibrium provides the following stability result
\cite{PZEJP13}:

\begin{proposition} \label{prop:Heterogeneous SIS stability condition}Consider
the mean-field approximation of the HeNeSIS model in (\ref{eq:HeNiSIS dynamics}) and assume that
$A_{\mathcal{G}}\geq0$, $\beta_{i},\delta_{i}>0$. Then, if the eigenvalue
with largest real part of $BA_{\mathcal{G}}-D$ satisfies 
\begin{equation}
\Re\left[\lambda_{1}\left(BA_{\mathcal{G}}-D\right)\right]\leq-\varepsilon,\label{eq:Spectral Control}
\end{equation}
for some $\varepsilon>0$, the disease-free equilibrium ($\boldsymbol{p}^{*}=\boldsymbol{0}$)
is globally exponentially stable, i.e., $\left\Vert \boldsymbol{p}\left(t\right)\right\Vert \leq\left\Vert \boldsymbol{p}\left(0\right)\right\Vert K\exp\left(-\varepsilon t\right)$,
for some $K>0$.

\end{proposition} 
\begin{rem}
In the proof of Proposition \ref{prop:Heterogeneous SIS stability condition}
in \cite{PZEJP13}, we show that the linear dynamical system $\dot{\boldsymbol{p}}\left(t\right)=\left(BA_{\mathcal{G}}-D\right)\boldsymbol{p}\left(t\right)$
upper-bounds the dynamics in (\ref{eq:HeNiSIS dynamics});
thus, the spectral result in (\ref{eq:Spectral Control}) is a sufficient
condition for the mean-field approximation of the HeNeSIS model to be globally exponentially
stable.
\end{rem}

\subsection{\label{sub:Problem-Statements}Problem Statements}

We describe two resource allocation problems to contain the
spread of an infection by distributing protection
resources throughout the network. We consider two types
of protection resources: (\emph{i}) preventive resources (or vaccinations),
and (\emph{ii}) corrective resources (or antidotes). Allocating preventive
resources at node $v_{i}$ reduces the infection rate $\beta_{i}$.
Allocating corrective resources at node $v_{i}$ increases the recovery
rate $\delta_{i}$. We assume that we are able to, simultaneously,
modify the infection and recovery rates of $v_{i}$ within feasible
intervals $0<\underline{\beta}_{i}\leq\beta_{i}\leq\bar{\beta}_{i}$
and $0<\underline{\delta}_{i}\leq\delta_{i}\leq\bar{\delta}_{i}$.
We consider that protection resources have an associated
cost. We define two cost functions, the vaccination cost function
$f_{i}\left(\beta_{i}\right)$ and the antidote cost function $g_{i}\left(\delta_{i}\right)$,
that account for the cost of tuning the infection and recovery rates
of node $v_{i}$ to $\beta_{i}\in\left[\underline{\beta}_{i},\bar{\beta}_{i}\right]$
and $\delta_{i}\in\left[\underline{\delta}_{i},\bar{\delta}_{i}\right]$,
respectively. In the rest of the paper we assume that the vaccination
cost function $f_{i}\left(\beta_{i}\right)$ is monotonically decreasing
w.r.t. $\beta_{i}$ and the antidote cost function $g_{i}\left(\delta_{i}\right)$
is monotonically increasing w.r.t. $\delta_{i}$.

In this context, we study two types of resource allocation problems
for the HeNiSIS model: (\emph{i}) the \emph{rate-constrained} allocation
problem, and (\emph{ii}) the \emph{budget-constrained} allocation
problem. In the \emph{rate-constrained} problem, we find the cost-optimal
distribution of vaccines and antidotes to achieve a given exponential
decay rate in the vector of infections, i.e., given $\overline{\varepsilon}$,
allocate resources such that $\left\Vert \boldsymbol{p}\left(t\right)\right\Vert \leq\left\Vert \boldsymbol{p}\left(0\right)\right\Vert K\exp\left(-\overline{\varepsilon}t\right)$,
$K>0$. In the\emph{ budget-constrained} problem, we are given a total
budget $C$ and we find the best allocation of vaccines and/or antidotes
to maximize the exponential decay rate of $\left\Vert \boldsymbol{p}\left(t\right)\right\Vert $,
i.e., maximize $\varepsilon$ (the decay rate) such that $\left\Vert \boldsymbol{p}\left(0\right)\right\Vert K\exp\left(-\overline{\varepsilon}t\right)$.

Based on Proposition \ref{prop:Heterogeneous SIS stability condition},
the decay rate of an epidemic outbreak is determined by $\varepsilon$
in (\ref{eq:Spectral Control}). Thus, we can formulate the rate-constrained
problem, as follows:
\begin{problem}
\emph{\label{Problem:Rate Constrained Allocation}(Rate-constrained}
\emph{allocation) Given the following elements: (i) A (positively)
weighted, directed network $\mathcal{G}$ with adjacency matrix $A_{\mathcal{G}}$,
(ii) a set of cost functions $\left\{ f_{i}\left(\beta_{i}\right),g_{i}\left(\delta_{i}\right)\right\} _{i=1}^{n}$,
(iii) bounds on the infection and recovery rates $0<\underline{\beta}_{i}\leq\beta_{i}\leq\overline{\beta}_{i}$
and $0<\underline{\delta}_{i}\leq\delta_{i}\leq\overline{\delta}_{i}$,
$i=1,\ldots,n$, and (iv) a desired exponential decay rate $\overline{\varepsilon}>0$;
find the cost-optimal distribution of vaccines and antidotes to achieve
the desired decay rate.}
\end{problem}
Given a desired decay rate $\overline{\varepsilon}$, the rate-constrained
allocation problem can be stated as the following optimization problem:

\emph{
\begin{align}
\underset{{\scriptscriptstyle \left\{ \beta_{i},\delta_{i}\right\} _{i=1}^{n}}}{\mbox{minimize }} & \sum_{i=1}^{n}f_{i}\left(\beta_{i}\right)+g_{i}\left(\delta_{i}\right)\label{eq:Rate-Constrain Spectral Problem}\\
\mbox{subject to } & \Re\left[\lambda_{1}\left(\mbox{diag}\left(\beta_{i}\right)A_{\mathcal{G}}-\mbox{diag}\left(\delta_{i}\right)\right)\right]\leq-\overline{\varepsilon},\label{eq:Spectral constraint}\\
 & \underline{\beta}_{i}\leq\beta_{i}\leq\overline{\beta}_{i},\label{eq:Square constraint beta}\\
 & \underline{\delta}_{i}\leq\delta_{i}\leq\overline{\delta}_{i},\mbox{ }i=1,\ldots,n,\label{eq:Square constraint delta}
\end{align}
}where (\ref{eq:Rate-Constrain Spectral Problem}) is the total investment,
(\ref{eq:Spectral constraint}) constrains the decay rate to $\overline{\varepsilon}$,
and (\ref{eq:Square constraint beta})-(\ref{eq:Square constraint delta})
maintain the infection and recovery rates in their feasible limits.

Similarly, given a budget $C$, the budget-constrained allocation
problem is formulated as follows:
\begin{problem}
\label{Problem: Budget Constrained Allocation}(\emph{Budget-constrained
allocation}) \emph{Given the following elements: (i) A (positively)
weighted, directed network $\mathcal{G}$ with adjacency matrix $A_{\mathcal{G}}$,
(ii) a set of cost functions $\left\{ f_{i}\left(\beta_{i}\right),g_{i}\left(\delta_{i}\right)\right\} _{i=1}^{n}$,
}(\emph{iii})\emph{ bounds on the infection and recovery rates $0<\underline{\beta}_{i}\leq\beta_{i}\leq\overline{\beta}_{i}$
and $0<\underline{\delta}_{i}\leq\delta_{i}\leq\overline{\delta}_{i}$,
$i=1,\ldots,n$,} \emph{and (iv) a total budget $C$; find the cost-optimal
distribution of vaccines and antidotes to maximize the exponential
decay rate $\varepsilon$.}
\end{problem}
Based on Proposition \ref{prop:Heterogeneous SIS stability condition},
we can state this problem as the following optimization program:

\emph{
\begin{align}
\underset{{\scriptscriptstyle \varepsilon,\left\{ \beta_{i},\delta_{i}\right\} _{i=1}^{n}}}{\mbox{maximize }} & \varepsilon\label{eq:Budget-Constrained Spectral Problem}\\
\mbox{subject to } & \Re\left[\lambda_{1}\left(\mbox{diag}\left(\beta_{i}\right)A_{\mathcal{G}}-\mbox{diag}\left(\delta_{i}\right)\right)\right]\leq-\varepsilon,\label{eq:Spectral constraint in budget problem}\\
 & \sum_{i=1}^{n}f_{i}\left(\beta_{i}\right)+g_{i}\left(\delta_{i}\right)\leq C,\label{eq:Budget constraint in budget problem}\\
 & \underline{\beta}_{i}\leq\beta_{i}\leq\overline{\beta}_{i},\label{eq:Square constraint for beta in budget problem}\\
 & \underline{\delta}_{i}\leq\delta_{i}\leq\overline{\delta}_{i},\mbox{ }i=1,\ldots,n,\label{eq:Square constraint for delta in budget problem}
\end{align}
}where (\ref{eq:Budget constraint in budget problem}) is the budget
constraint.

In the following section, we propose an approach to solve these problems
in polynomial time for weighted and directed contact networks, under
certain assumptions on the cost functions $f_{i}$ and $g_{i}$.

\section{\label{sec:Convex Framework}A Convex Framework for Optimal Resource
Allocation}

We propose a convex formulation to solve both the budget-constrained
and the rate-constrained allocation problem in weighted, directed
networks using \emph{geometric programming (GP)} \cite{BV04}. We
first provide a solution for strongly connected digraphs in Subsection
\ref{sub:GP for SC Digraphs}. We then extend our results to general
digraphs (not necessarily strongly connected) in Subsection \ref{sub:Allocation-for general digraphs}.

We start our exposition by briefly reviewing some
concepts used in our formulation. Let $x_{1},\ldots,x_{n}>0$ denote
$n$ decision variables and define $\mathbf{x}\triangleq\left(x_{1},\ldots,x_{n}\right)\in\mathbb{R}_{++}^{n}$.
In the context of GP, a \emph{monomial $h(\mathbf{x})$} is defined
as a real-valued function of the form $h(\mathbf{x})\triangleq dx_{1}^{a_{1}}x_{2}^{a_{2}}\ldots x_{n}^{a_{n}}$
with $d>0$ and $a_{i}\in\mathbb{R}$. A \emph{posynomial} function
$q(\mathbf{x})$ is defined as a sum of monomials, i.e., $q(\mathbf{x})\triangleq\sum_{k=1}^{K}c_{k}x_{1}^{a_{1k}}x_{2}^{a_{2k}}\ldots x_{n}^{a_{nk}}$,
where $c_{k}>0$.

In our formulation, it is useful to define the following class of
functions:
\begin{defn}
A function $F:\mathbb{R}^{n}\to\mathbb{R}$ is \emph{convex in log-scale}
if the function 
\begin{equation}
F\left(\mathbf{y}\right)\triangleq\log f\left(\exp\mathbf{y}\right),\label{eq:Convex in Log Scale}
\end{equation}
is convex in \textbf{$\mathbf{y}$} (where $\exp\mathbf{y}$ indicates
component-wise exponentiation).\end{defn}
\begin{rem}
Note that posynomials (hence, also monomials) are convex in log-scale
\cite{BV04}.
\end{rem}
A geometric program (GP) is an optimization problem of the form (see
\cite{BKVH07} for a comprehensive treatment):
\begin{align}
\mbox{minimize } & f(\mathbf{x})\label{eq:General GP}\\
\mbox{subject to } & q_{i}(\mathbf{x})\leq1,\: i=1,...,m,\nonumber \\
 & h_{i}(\mathbf{x})=1,\: i=1,...,p,\nonumber 
\end{align}
where $q_{i}$ are posynomial functions, $h_{i}$ are monomials, and
$f$ is a convex function in log-scale%
\footnote{Geometric programs in standard form are usually formulated assuming
$f\left(\mathbf{x}\right)$ is a posynomial. In our formulation, we
assume that $f\left(\mathbf{x}\right)$ is in the broader class of
convex functions in logarithmic scale.%
}. A GP is a quasiconvex optimization problem \cite{BV04} that can
be transformed to a convex problem. This conversion is based on
the logarithmic change of variables $y_{i}=\log x_{i}$, and a logarithmic
transformation of the objective and constraint functions (see \cite{BKVH07}
for details on this transformation). After this transformation, the
GP in (\ref{eq:General GP}) takes the form
\begin{align}
\mbox{minimize } & F\left(\mathbf{y}\right)\label{eq:Transformed GP}\\
\mbox{subject to } & Q_{i}\left(\mathbf{y}\right)\leq0,\: i=1,...,m,\nonumber \\
 & \mathbf{b}_{i}^{T}\mathbf{y}+\log d_{i}=0,\: i=1,...,p,\nonumber 
\end{align}
where $Q_{i}\left(\mathbf{y}\right)\triangleq\log q_{i}(\exp\mathbf{y})$ and
$F\left(\mathbf{y}\right)\triangleq\log f\left(\exp\mathbf{y}\right)$. Also, assuming that $h_{i}\left(\mathbf{x}\right)\triangleq d_{i}x_{1}^{b_{1,i}}x_{2}^{b_{2,i}}\ldots x_{n}^{b_{n,i}}$,
we obtain the equality constraint above, with $\mathbf{b}_{i}\triangleq\left(b_{1,i},\ldots,b_{n,i}\right)$,
after the logarithmic change of variables. Notice
that, since $f\left(\mathbf{x}\right)$ is convex in log-scale, $F\left(\mathbf{y}\right)$
is a convex function. Also, since $q_{i}$ is a posynomial (therefore,
convex in log-scale), $Q_{i}$ is also a convex function. In conclusion,
(\ref{eq:Transformed GP}) is a convex optimization problem in standard
form and can be efficiently solved in polynomial time \cite{BV04}.

As we shall show in Subsections \ref{sub:GP for SC Digraphs} and
\ref{sub:Allocation-for general digraphs}, we can solve Problems
\ref{Problem:Rate Constrained Allocation} and \ref{Problem: Budget Constrained Allocation}
using GP if the cost function $\sum_{i=1}^{n}f_{i}\left(\beta_{i}\right)+g_{i}\left(\delta_{i}\right)$
is convex in log-scale. In practical applications, we model the individual
cost functions $f_{i}\left(\beta_{i}\right)$ and $g_{i}\left(\delta_{i}\right)$
as posynomials. In practice, posynomials functions can be used to fit any function that is convex in log-log scale with arbitrary accuracy. Furthermore, there are well-developed numerical methods to fit posynomials to real data (see \cite{BKVH07}, Section 8, for a treatment about
the modeling abilities of monomials and posynomials).

In the following sections, we show how to transform Problems \ref{Problem:Rate Constrained Allocation}
and \ref{Problem: Budget Constrained Allocation} into GP's. In our
transformation, we use the theory of nonnegative matrices and the
Perron-Frobenius lemma. Our derivations are different if the contact
graph $\mathcal{G}$ is strongly connected or not. We cover the case
of $\mathcal{G}$ being a strongly connected digraph in Subsection
\ref{sub:GP for SC Digraphs} and we extend the theory to general
digraphs in Subsection \ref{sub:Allocation-for general digraphs}.

\subsection{GP for Strongly Connected Digraphs\label{sub:GP for SC Digraphs}}

In our derivations, we use Perron-Frobenius lemma, from the theory
of nonnegative matrices \cite{meyer2000matrix}:
\begin{lem}
\label{lem:Perron-Frobenius}(Perron-Frobenius) Let $M$ be a nonnegative,
irreducible matrix. Then, the following statements about its spectral
radius, $\rho\left(M\right)$, hold:

\emph{(}a\emph{)} $\rho\left(M\right)>0$ is a simple eigenvalue of
$M$,

\emph{(}b\emph{)} $M\mathbf{u}=\rho\left(M\right)\mathbf{u}$, for
some $\mathbf{u}\in\mathbb{R}_{++}^{n}$, and

\emph{(}c\emph{)} \textup{$\rho\left(M\right)=\inf\left\{ \lambda\in\mathbb{R}:M\mathbf{u}\leq\lambda\mathbf{u}\mbox{ for }\mathbf{u}\in\mathbb{R}_{++}^{n}\right\} $.}\end{lem}
\begin{rem}
Since a matrix $M$ is \emph{irreducible} if and only if its associated
digraph $\mathcal{G}_{M}$ is strongly connected, the above lemma
also holds for the spectral radius of the adjacency matrix of any
(positively) weighted, strongly connected digraph.
\end{rem}
From Lemma \ref{lem:Perron-Frobenius}, we infer the following results:
\begin{cor}
\label{cor:Eig equals Rad}Let $M$ be a nonnegative, irreducible
matrix. Then, its eigenvalue with the largest real part, $\lambda_{1}\left(M\right)$,
is real, simple, and equal to the spectral radius $\rho\left(M\right)>0$.\end{cor}
\begin{lem}
\label{lem:Monotonicity of lambda1}Consider the adjacency matrix
$A_{\mathcal{G}}$ of a (positively) weighted, directed, strongly
connected graph $\mathcal{G}$, and two sets of positive numbers $\left\{ \beta_{i}\right\} _{i=1}^{n}$
and $\left\{ \delta_{i}\right\} _{i=1}^{n}$. Then, $\lambda_{1}\left(\mbox{diag}\left(\beta_{i}\right)A-\mbox{diag}\left(\delta_{i}\right)\right)$
is an increasing function w.r.t. $\beta_{k}$ (respectively, monotonically
decreasing w.r.t. $\delta_{k}$) for $k=1,\ldots,n$.\end{lem}
\begin{IEEEproof}
In the Appendix.
\end{IEEEproof}
From the above results, we have the following result (\cite{BV04},
Chapter 4):
\begin{prop}
\label{prop:From PF to Posynomials}Consider the $n\times n$ nonnegative,
irreducible matrix $M\left(\mathbf{x}\right)$ with entries being
either $0$ or posynomials with domain $\mathbf{x}\in\mathcal{S}\subseteq\mathbb{R}_{++}^{k}$,
where $\mathcal{S}$ is defined as $\mathcal{S}=\bigcap_{i=1}^{m}\left\{ \mathbf{x\in\mathbb{R}}_{++}^{k}:f_{i}\left(\mathbf{x}\right)\leq1\right\} $,
$f_{i}$ being posynomials. Then, we can minimize $\lambda_{1}\left(M\left(\mathbf{x}\right)\right)$
for $\mathbf{x}\in\mathcal{S}$ solving the following GP:
\begin{align}
\underset{\lambda,\left\{ u_{i}\right\} _{i=1}^{n},\mathbf{x}}{\mbox{minimize }} & \lambda\label{eq:GP for spectral objective}\\
\mbox{subject to } & \frac{\sum_{j=1}^{n}M_{ij}\left(\mathbf{x}\right)u_{j}}{\lambda u_{i}}\leq1,\: i=1,\ldots,n,\label{eq:GP for spectral constraint}\\
 & f_{i}\left(\mathbf{x}\right)\leq1,\: i=1,\ldots,m.
\end{align}

\end{prop}

Based on the above results, we provide solutions to both the rate-constrained
and the budget-constrained problems.

\subsubsection{\label{sub:Solution to Budget Constrained}Solution to the Budget-Constrained
Allocation Problem for Strongly Connected Digraphs}

Assuming that the cost functions $f_{i}$ and $g_{i}$ are posynomials
and the contact graph $\mathcal{G}$ is strongly connected, the budget-constrained
allocation problem in \ref{Problem: Budget Constrained Allocation}
can be solved as follows:
\begin{thm}
\label{thm:GP for budget constrained}Consider the following elements:
(i) A strongly connected graph $\mathcal{G}$ with adjacency matrix
$A_{\mathcal{G}}=[A_{ij}]$, (ii) posynomial cost functions \emph{$\left\{ f_{i}\left(\beta_{i}\right),g_{i}\left(\delta_{i}\right)\right\} _{i=1}^{n}$},
(iii) bounds on the infection and recovery rates $0<\underline{\beta}_{i}\leq\beta_{i}\leq\overline{\beta}_{i}$
and \textup{$0<\underline{\delta}_{i}\leq\delta_{i}\leq\overline{\delta}_{i}$},
$i=1,\ldots,n$, and (iv) a maximum budget $C$ to invest in protection
resources. Then, the optimal investment on vaccines and antidotes
for node $v_{i}$ to solve Problem \ref{Problem: Budget Constrained Allocation}
are $f_{i}\left(\beta_{i}^{*}\right)$ and $g_{i}\left(\overline{\Delta}+1-\widehat{\delta}_{i}^{*}\right)$,
where $\overline{\Delta}\triangleq\max\left\{ \overline{\delta}_{i}\right\} _{i=1}^{n}$
and \textup{\emph{$\beta_{i}^{*}$,$\widehat{\delta}_{i}^{*}$ are
the optimal solution for $\beta_{i}$ and $\widehat{\delta}_{i}$
in the following GP}}:\emph{
\begin{align}
\underset{{\scriptstyle \lambda,\left\{ u_{i},\beta_{i},\widehat{\delta}_{i},t_{i}\right\} _{i=1}^{n}}}{\mbox{minimize}} & \lambda\label{eq:Budget-Constrained Spectral Problem-1}\\
\mbox{subject to } & \frac{\beta_{i}\sum_{j=1}^{n}A_{ij}u_{j}+\widehat{\delta}_{i}u_{i}}{\lambda u_{i}}\leq1,\label{eq:Eigenvalu condition in spectral constraint}\\
 & \sum_{k=1}^{n}f_{k}\left(\beta_{k}\right)+g_{k}\left(t_{k}\right)\leq C,\label{eq:budget constraint in budget constrained}\\
 & \left(t_{i}+\widehat{\delta}_{i}\right)\left/\left(\overline{\Delta}+1\right)\right.\leq1,\label{eq:t trick in budget constrained}\\
 & \overline{\Delta}+1-\overline{\delta}_{i}\leq\widehat{\delta}_{i}\leq\overline{\Delta}+1-\underline{\delta}_{i},\label{eq:Beta contraint in budget constrained}\\
 & \underline{\beta}_{i}\leq\beta_{i}\leq\overline{\beta}_{i},\mbox{ }i=1,\ldots,n.\label{eq:Delta constraint in budget constraint}
\end{align}
}\end{thm}
\begin{IEEEproof}
First, based on Proposition \ref{prop:From PF to Posynomials}, we
have that maximizing $\varepsilon$ in (\ref{eq:Budget-Constrained Spectral Problem})
subject to (\ref{eq:Spectral constraint in budget problem})-(\ref{eq:Square constraint for beta in budget problem})
is equivalent to minimizing $\lambda_{1}\left(BA_{\mathcal{G}}-D\right)$
subject to (\ref{eq:Budget constraint in budget problem}) and (\ref{eq:Square constraint for beta in budget problem}),
where $B\triangleq\mbox{diag}\left(\beta_{i}\right)$ and $D\triangleq\mbox{diag}\left(\delta_{i}\right)$.
Let us define $\widehat{D}\triangleq\mbox{diag}\left(\widehat{\delta}_{i}\right)$,
where $\widehat{\delta}_{i}\triangleq\overline{\Delta}+1-\delta_{i}$
and $\overline{\Delta}\triangleq\max\left\{ \overline{\delta}_{i}\right\} _{i=1}^{n}$.
Then, $\lambda_{1}\left(BA_{\mathcal{G}}+\widehat{D}\right)=\lambda_{1}\left(BA_{\mathcal{G}}-D\right)+\overline{\Delta}+1$.
Hence, minimizing $\lambda_{1}\left(BA_{\mathcal{G}}-D\right)$ is
equivalent to minimizing $\lambda_{1}\left(BA_{\mathcal{G}}+\widehat{D}\right)$.
The matrix $BA_{\mathcal{G}}+\widehat{D}$ is nonnegative and irreducible
if $A_{\mathcal{G}}$ is the adjacency matrix of a strongly connected
digraph. Therefore, applying Proposition \ref{prop:From PF to Posynomials},
we can minimize $\lambda_{1}\left(BA_{\mathcal{G}}+\widehat{D}\right)$
by minimizing the cost function in (\ref{eq:Budget-Constrained Spectral Problem-1})
under the constraints (\ref{eq:Eigenvalu condition in spectral constraint})-(\ref{eq:Delta constraint in budget constraint}).
Constraints (\ref{eq:Beta contraint in budget constrained}) and (\ref{eq:Delta constraint in budget constraint})
represent bounds on the achievable infection and curing rates. Notice
that we also have a constraint associated to the budget available,
i.e., $\sum_{k=1}^{n}f_{k}\left(\beta_{k}\right)+g_{k}\left(\overline{\Delta}+1-\widehat{\delta}_{i}\right)\leq C$.
But, even though $g_{k}$$\left(\delta_{k}\right)$ is a polynomial
function in $\delta_{k}$, $g_{k}\left(\overline{\Delta}+1-\widehat{\delta}_{k}\right)$
is not a posynomial in $\widehat{\delta}_{i}$. To overcome this issue,
we can replace the argument of $g_{k}$ by a new variable $t_{k}$,
along with the constraint $t_{k}\leq\overline{\Delta}+1-\widehat{\delta}_{k}$,
which can be expressed as the posynomial inequality, $\left(t_{k}+\widehat{\delta}_{k}\right)/\left(\overline{\Delta}+1\right)\leq1$.
As we proved in Lemma \ref{lem:Monotonicity of lambda1}, the largest
eigenvalue $\lambda_{1}\left(BA-D\right)$ is a decreasing value of
$\delta_{k}$ and the antidote cost function $g_{k}$ is monotonically
increasing w.r.t. $\delta_{k}$. Thus, adding the inequality $t_{k}\leq\overline{\Delta}+1-\widehat{\delta}_{k}$
does not change the result of the optimization problem, since at optimality
$t_{k}$ will saturate to its largest possible value $t_{k}=\overline{\Delta}+1-\widehat{\delta}_{k}$.
\end{IEEEproof}

\subsubsection{\label{sub:Solution to Rate Constrained}Solution to Rate-Constrained
Allocation Problem for Strongly Connected Digraphs}

Problem \ref{Problem:Rate Constrained Allocation} can be written
as the following optimization program:
\begin{thm}
\label{thm:GP for rate constrained}Consider the following elements:
(i) A strongly connected graph $\mathcal{G}$ with adjacency matrix
$A_{\mathcal{G}}=[A_{ij}]$, (ii) posynomial cost functions \emph{$\left\{ f_{i}\left(\beta_{i}\right),g_{i}\left(\delta_{i}\right)\right\} _{i=1}^{n}$},
(iii) bounds on the infection and recovery rates $0<\underline{\beta}_{i}\leq\beta_{i}\leq\overline{\beta}_{i}$
and \textup{$0<\underline{\delta}_{i}\leq\delta_{i}\leq\overline{\delta}_{i}$},
$i=1,\ldots,n$, and (iv) a desired exponential decay rate $\overline{\varepsilon}$.
Then, the optimal investment on vaccines and antidotes for node\emph{
$v_{i}$ }to solve Problem \ref{Problem:Rate Constrained Allocation}
are $f_{i}\left(\beta_{i}^{*}\right)$ and $g_{i}\left(\widetilde{\Delta}+1-\widetilde{\delta}_{i}^{*}\right)$,
where $\widetilde{\Delta}\triangleq\max\left\{ \overline{\varepsilon},\overline{\delta}_{i}\mbox{ for }i=1,\ldots,n\right\} $
and \textup{\emph{$\beta_{i}^{*}$,}}\textup{$\widetilde{\delta}_{i}^{*}$}\textup{\emph{
are the optimal solution for $\beta_{i}$ and $\widetilde{\delta}_{i}$
in the following GP}}:

\emph{
\begin{align}
\underset{{\scriptstyle \left\{ u_{i},\beta_{i},\widetilde{\delta}_{i},t_{i}\right\} _{i=1}^{n}}}{\mbox{minimize}} & \sum_{k=1}^{n}f_{k}\left(\beta_{k}\right)+g_{k}\left(t_{k}\right)\label{eq:Budget-Constrained Spectral Problem-1-1}\\
\mbox{subject to } & \frac{\beta_{i}\sum_{j=1}^{n}A_{ij}u_{j}+\widetilde{\delta}_{i}u_{i}}{\left(\widetilde{\Delta}+1-\overline{\varepsilon}\right)u_{i}}\leq1,\label{eq:Spectral constraint in budget constrained}\\
 & \left(t_{i}+\widetilde{\delta}_{i}\right)\left/\left(\widetilde{\Delta}+1\right)\right.\leq1,\label{eq:t trick in rate constrained}\\
 & \widetilde{\Delta}+1-\overline{\delta}_{i}\leq\widehat{\delta}_{i}\leq\widetilde{\Delta}+1-\underline{\delta}_{i},\label{eq:delta constraint in rate constrained}\\
 & \underline{\beta}_{i}\leq\beta_{i}\leq\overline{\beta}_{i},\, i=1,\ldots,n.\label{eq:beta constraint in rate constrained}
\end{align}
}\end{thm}
\begin{IEEEproof}
(The proof is similar to the one for Theorem \ref{thm:GP for budget constrained}
and we only present here the main differences.) Define $\widetilde{D}\triangleq\mbox{diag}\left(\widetilde{\delta}_{i}\right)$
where $\widetilde{\delta}_{i}\triangleq\widetilde{\Delta}+1-\delta_{i}$
and $\widetilde{\Delta}\triangleq\max\left\{ \overline{\varepsilon},\overline{\delta}_{i}\mbox{ for }i=1,\ldots,n\right\} $.
Since $\lambda_{1}\left(BA_{\mathcal{G}}+\widetilde{D}\right)=\lambda_{1}\left(BA_{\mathcal{G}}-D\right)+\widetilde{\Delta}+1$,
the spectral condition $\lambda_{1}\left(BA_{\mathcal{G}}-D\right)\leq-\overline{\varepsilon}$
is equivalent to $\lambda_{1}\left(BA_{\mathcal{G}}+\widetilde{D}\right)\leq\widetilde{\Delta}+1-\overline{\varepsilon}$.
From the definition of $\widetilde{\Delta}$ we have that $\widetilde{\Delta}+1-\overline{\varepsilon}>0$.
Also, $BA_{\mathcal{G}}+\widetilde{D}$ is a nonnegative and irreducible
matrix if $\mathcal{G}$ is a strongly connected digraph. From (\ref{eq:GP for spectral constraint}),
we can write the spectral constraint $\lambda_{1}\left(BA_{\mathcal{G}}+\widetilde{D}\right)\leq\widetilde{\Delta}+1-\overline{\varepsilon}$
as
\[
\frac{\beta_{i}\sum_{j=1}^{n}A_{ij}u_{j}+\widetilde{\delta}_{i}u_{i}}{\left(\widetilde{\Delta}+1-\overline{\varepsilon}\right)u_{i}}\leq1,
\]
for $u_{i}\in\mathbb{R}_{++}$, $\lambda\in\mathbb{R}$, which results
in constraint (\ref{eq:Spectral constraint in budget constrained}).
The rest of constraints can be derived following similar derivations
as in the Proof of Theorem \ref{thm:GP for budget constrained}.
\end{IEEEproof}
In Subsections \ref{sub:Solution to Budget Constrained} and \ref{sub:Solution to Rate Constrained},
we have presented two geometric programs to find the optimal solutions
to both the budget-constrained and the rate-constrained allocation
problems. In our derivations, we have made the assumption of $\mathcal{G}$
being a strongly connected graph. In the next subsection, we show
how to solve these allocation problems for any digraphs, after relaxing
the strong connectivity assumption.

\subsection{\label{sub:Allocation-for general digraphs}Solution to Allocation
Problems for General Digraphs}

The Perron-Frobenius lemma state that given a nonnegative, \emph{irreducible}
matrix $M$, its spectral radius $\rho(M)$ is simple and strictly
positive (thus, $\rho\left(M\right)=\lambda_{1}\left(M\right)$) and
the associated dominant eigenvector has strictly positive components.
Perron-Frobenius lemma is not applicable to digraphs that are not
strongly connected, since the associated adjacency matrix is not irreducible.
For weighted (possibly reducible) digraphs, the statements in the
Perron-Frobenius lemma are weaken, as follows \cite{meyer2000matrix}:
\begin{lem}
Let $M$ be a nonnegative matrix. Then, the following statements hold:\end{lem}
\begin{description}
\item [{\textmd{(}\textmd{\emph{a}}\textmd{)}}] $\rho\left(M\right)\geq0$
is an eigenvalue of $M$ (not necessarily simple).
\item [{\textmd{(}\textmd{\emph{b}}\textmd{)}}] $M\mathbf{u}=\rho\left(M\right)\mathbf{u}$,
for some $\mathbf{u}\in\mathbb{R}_{+}^{n}$.
\item [{\textmd{(}\textmd{\emph{c}}\textmd{)}}] $\rho\left(M\right)=\inf\left\{ \lambda\in\mathbb{R}:M\mathbf{u}\leq\lambda\mathbf{u}\mbox{ for }\mathbf{u}\in\mathbb{R}_{+}^{n}\right\} $.\end{description}
\begin{rem}
Notice that in item (\emph{c}), the components of $\mathbf{u}$ are
nonnegative (instead of positive). This is an issue if we want to
use Proposition \ref{prop:From PF to Posynomials}, since the components
of $\mathbf{v}$ must be strictly positive to use GP. In what follows,
we show how to resolve this issue.
\end{rem}
Let us define the function $\mathcal{Z}\left(\mathbf{u}\right)\triangleq\left\{ i:u_{i}=0\right\} $,
i.e., a function that returns the set of indexes indicating the location
of the zero entries of a vector $\mathbf{u}=\left[u_{i}\right]$.
\begin{lem}
\label{lem:Pattern of Zeros}Consider a square matrix $M$. The following
transformations preserve the location of zeros in the dominant eigenvector:

(a) $T_{\alpha}:M\to M+\alpha I$, for any $\alpha\in\mathbb{R}$,
and

(b) $T_{R}:M\to RM$, for $M\geq0$ and $R=\mbox{diag}\left(r_{i}\right)$,
$r_{i}>0$.\end{lem}
\begin{IEEEproof}
In the Appendix.\end{IEEEproof}
\begin{prop}
\label{prop:Zeros in BA-D}Consider a nonnegative matrix $A$ and
two diagonal matrices $B=\mbox{diag}\left(b_{i}\right)$ and $D=\mbox{diag}\left(d_{i}\right)$
with $b_{i},d_{i}>0$. Then, the location of the zero entries of the
dominant eigenvector of $BA-D$ are the same as those of $A$, i.e.,
$\mathcal{Z}\left(\mathbf{v}_{1}\left(BA-D\right)\right)=\mathcal{Z}\left(\mathbf{v}_{1}\left(A\right)\right)$.\end{prop}
\begin{IEEEproof}
In the Appendix.
\end{IEEEproof}
Proposition \ref{prop:Zeros in BA-D} allows us to know the location
of the zeros of $\mathbf{v}_{1}\left(BA_{\mathcal{G}}-D\right)$ for
any given graph $A_{\mathcal{G}}\geq0$, independently of the allocation
of vaccines and antidotes in the network. This location is determined
by the set $\mathcal{Z}\left(\mathbf{v}_{1}\left(A_{\mathcal{G}}\right)\right)\triangleq\mathcal{Z}_{\mathcal{G}}$,
which is the set of nodes with zero eigenvector centrality \cite{New10}.
Hence, we can exclude the variables $u_{i}$ for $i\in\mathcal{Z}_{\mathcal{G}}$
from the GP's in Theorems \ref{thm:GP for rate constrained} and \ref{thm:GP for budget constrained}.
Hence, since the components in the set $\left\{ u_{i}:i\in\mathcal{Z}_{\mathcal{G}}\right\} $
are not part of the spectral conditions (\ref{eq:Eigenvalu condition in spectral constraint})
and (\ref{eq:Spectral constraint in budget constrained}), we can
split the allocation problems into two different sets of decision
variables. We elaborate on this splitting in the following subsections.

\subsubsection{\label{sub:Rate-Constrained-for general digraphs}Rate-Constrained
Allocation Problem for General Digraphs}

The set of decision variables in (\ref{eq:Budget-Constrained Spectral Problem-1-1})
split into two sets: $V_{z}\triangleq\{u_{i},\beta_{i},\widetilde{\delta}_{i},t_{i}\}_{i\in\mathcal{Z_{G}}}$
and $V_{nz}\triangleq\{u_{i},\beta_{i},\widetilde{\delta}_{i},t_{i}\}_{i\notin\mathcal{Z_{G}}}$.
From (\ref{eq:Budget-Constrained Spectral Problem-1-1}), the following
optimization problem holds for the variables in $V_{z}$: 
\begin{align*}
\underset{{\scriptstyle \left\{ \beta_{i},\widetilde{\delta}_{i},t_{i}\right\} _{i\in\mathcal{Z}_{\mathcal{G}}}}}{\mbox{minimize}} & \sum_{k=1}^{n}f_{k}\left(\beta_{k}\right)+g_{k}\left(\widetilde{\Delta}+1-\widetilde{\delta}_{i}^{*}\right)\\
\mbox{subject to } & \widetilde{\Delta}+1-\overline{\delta}_{i}\leq\widehat{\delta}_{i}\leq\widetilde{\Delta}+1-\underline{\delta}_{i},\\
 & \underline{\beta}_{i}\leq\beta_{i}\leq\overline{\beta}_{i},\mbox{ for }i\in\mathcal{Z}_{\mathcal{G}}.
\end{align*}
Thus, for $f_{i}$ decreasing and $g_{i}$ increasing, it is easy
to verify that the optimal infection and recovery rates are $\beta_{i}^{*}=\overline{\beta}_{i}$
and $\delta_{i}^{*}=\underline{\delta}_{i}$ for all $i\in\mathcal{Z}_{\mathcal{G}}$.
Those rates correspond to the minimum possible value of investment
on those nodes. In other words, nodes with zero \emph{eigenvector
centrality \cite{New10}} receive the minimum possible value of investment.

On the other hand, for those decision variables in $V_{nz}$, we can
to adapt the GP formulation in Theorem \ref{thm:GP for rate constrained},
as indicated in the following Theorem.
\begin{thm}
\label{thm:General Solution for Rate Constrained-1}Consider the following
elements: (i) A positively weighted digraph with adjacency matrix
$A_{\mathcal{G}}$, (ii) posynomial cost functions \emph{$\left\{ f_{i}\left(\beta_{i}\right),g_{i}\left(\delta_{i}\right)\right\} _{i=1}^{n}$},
(iii) bounds on the infection and recovery rates $0<\underline{\beta}_{i}\leq\beta_{i}\leq\overline{\beta}_{i}$
and \textup{$0<\underline{\delta}_{i}\leq\delta_{i}\leq\overline{\delta}_{i}$},
$i=1,\ldots,n$, and (iv) a desired exponential decay rate $\overline{\varepsilon}$.
Then, the optimal spreading and recovery rate in Problem \ref{Problem:Rate Constrained Allocation}
are $\beta_{i}^{*}=\overline{\beta}_{i}$ and $\delta_{i}^{*}=\underline{\delta}_{i}$
for \textup{$i\in\mathcal{Z}_{\mathcal{G}}$. For $i\mathcal{\notin Z_{G}}$,
}the optimal rates can be computed from the optimal solution of the
following GP:

\emph{
\begin{align}
\underset{{\scriptstyle \left\{ u_{i},\beta_{i},\widetilde{\delta}_{i},t_{i}\right\} _{i\mathcal{\notin Z}_{\mathcal{G}}}}}{\mbox{minimize}} & \sum_{k=1}^{n}f_{k}\left(\beta_{k}\right)+g_{k}\left(t_{k}\right)\label{eq:Budget-Constrained Spectral Problem-1-1-1}\\
\mbox{subject to } & \frac{\beta_{i}\sum_{j\mathcal{\notin Z_{G}}}A_{ij}u_{j}+\widetilde{\delta}_{i}u_{i}}{\left(\widetilde{\Delta}+1-\overline{\varepsilon}\right)u_{i}}\leq1,\label{eq:Spectral constraint in budget constrained-1}\\
 & \left(t_{i}+\widetilde{\delta}_{i}\right)\left/\left(\widetilde{\Delta}+1\right)\right.\leq1,\\
 & \widetilde{\Delta}+1-\overline{\delta}_{i}\leq\widehat{\delta}_{i}\leq\widetilde{\Delta}+1-\underline{\delta}_{i},\\
 & \underline{\beta}_{i}\leq\beta_{i}\leq\overline{\beta}_{i},\mbox{ for }i\mathcal{\notin Z}_{\mathcal{G}}.
\end{align}
}The optimal spreading rate $\beta_{i}^{*}$ is directly obtained
from the solution, and the recovery rate is \emph{$\delta_{i}^{*}=\widetilde{\Delta}+1-\widetilde{\delta}_{i}^{*}$,
where $\widetilde{\Delta}\triangleq\max\left\{ \overline{\varepsilon},\overline{\delta}_{i}\mbox{ for }i=1,\ldots,n\right\} $.}\end{thm}
\begin{rem}
Since, $u_{i}=0$ if and only if $i\in\mathcal{Z_{G}}$, all the decision
variables in the above GP are strictly positive.
\end{rem}

\subsubsection{\label{sub:Budget-Constrained-for general digraphs}Budget-Constrained
Allocation Problem for General Digraphs}

In this case, one can also use the splitting techniques in Subsection
\ref{sub:Rate-Constrained-for general digraphs} to show that for
$i\in\mathcal{Z}_{\mathcal{G}}$ the optimal spreading and recovery
rates are $\beta_{i}=\overline{\beta}_{i}$ and $\delta_{i}=\underline{\delta}_{i}$.
This again corresponds to the minimum possible level of investment
for nodes with zero eigenvector centrality. We therefore allocate
a total amount equal to $f_{i}\left(\overline{\beta}_{i}\right)+g_{i}\left(\underline{\delta}_{i}\right)$
on each one of the nodes with zero eigenvalue centrality. As a result,
we should allocate the remaining budget
\begin{equation}
C-\sum_{i\in\mathcal{Z_{G}}}f_{i}\left(\overline{\beta}_{i}\right)+g_{i}\left(\underline{\delta}_{i}\right)\triangleq\overline{C}\label{eq:C-bar}
\end{equation}
to the set of nodes $\left\{ v_{i}\in\mathcal{V}:i\notin\mathcal{Z}_{\mathcal{G}}\right\} $.
Thus, the budget-constrained allocation problem in \ref{Problem: Budget Constrained Allocation}
can be written as the following GP for general digraphs:
\begin{thm}
\label{thm:General Solution for Budget Constrained}Consider the following
elements: (i) A positively weighted digraph with adjacency matrix
$A_{\mathcal{G}}$, (ii) posynomial cost functions \emph{$\left\{ f_{i}\left(\beta_{i}\right),g_{i}\left(\delta_{i}\right)\right\} _{i=1}^{n}$},
(iii) bounds on the infection and recovery rates $0<\underline{\beta}_{i}\leq\beta_{i}\leq\overline{\beta}_{i}$
and \textup{$0<\underline{\delta}_{i}\leq\delta_{i}\leq\overline{\delta}_{i}$},
$i=1,\ldots,n$, and (iv) a maximum budget $C$ to invest in protection
resources. Then, the optimal spreading and recovery rate in Problem
\ref{Problem: Budget Constrained Allocation} are $\beta_{i}^{*}=\overline{\beta}_{i}$
and $\delta_{i}^{*}=\underline{\delta}_{i}$ for \textup{$i\in\mathcal{Z}_{\mathcal{G}}$.
For $i\mathcal{\notin Z_{G}}$, }the optimal rates can be computed
from the optimal solution of the following GP:

\emph{
\begin{align}
\underset{{\scriptstyle \lambda,\left\{ u_{i},\beta_{i},\widehat{\delta}_{i},t_{i}\right\} _{i\mathcal{\notin Z_{G}}}}}{\mbox{minimize}} & \lambda\label{eq:Budget-Constrained Spectral Problem-1-2}\\
\mbox{subject to } & \frac{\beta_{i}\sum_{j\notin Z_{G}}A_{ij}u_{j}+\widehat{\delta}_{i}u_{i}}{\lambda u_{i}}\leq1,\label{eq:Eigenvalu condition in spectral constraint-1}\\
 & \sum_{k\mathcal{\notin Z_{G}}}f_{k}\left(\beta_{k}\right)+g_{k}\left(t_{k}\right)\leq\overline{C},\\
 & \left(t_{i}+\widehat{\delta}_{i}\right)\left/\left(\overline{\Delta}+1\right)\right.\leq1,\\
 & \overline{\Delta}+1-\overline{\delta}_{i}\leq\widehat{\delta}_{i}\leq\overline{\Delta}+1-\underline{\delta}_{i},\\
 & \underline{\beta}_{i}\leq\beta_{i}\leq\overline{\beta}_{i},\mbox{ }i\mathcal{\notin Z_{G}},
\end{align}
}where $\overline{C}$ is defined in (\ref{eq:C-bar}), the optimal
spreading rate $\beta_{i}^{*}$ is directly obtained from the solution,
and the recovery rate is \emph{$\delta_{i}^{*}=\overline{\Delta}+1-\widehat{\delta}_{i}^{*}$,
where $\overline{\Delta}\triangleq\max\left\{ \overline{\delta}_{i}\right\} _{i=1}^{n}$.}
\end{thm}
Theorems \ref{thm:General Solution for Rate Constrained-1}
and \ref{thm:General Solution for Budget Constrained} solve the optimal resource
allocation problems herein described for \emph{weighted, directed}
networks of \emph{nonidentical} agents.

\section{\label{sub:Numerical-Results}Numerical Results}

We apply our results to the design of a cost-optimal protection strategy
against epidemic outbreaks that propagate through the air transportation
network \cite{schneider2011suppressing}. We analyze real data from
the world-wide air transportation network and find the optimal distribution
of vaccines and antidotes to control (or contain) the spread of an epidemic outbreak. We consider both the rate-constraint and the budget-constrained
problems in our simulations. We limit our analysis to an air transportation
network spanning the major airports in the world; in particular, we
consider only airports having an incoming traffic greater than 10
million passengers per year (MPPY). There are $56$ such airports
world-wide and they are connected via $1,843$ direct flights, which
we represent as directed edges in a graph. To each directed edge,
we assign a weight equal to the number of passengers taking that flight
throughout the year (in MPPY units).

The weighted, directed graph representing the air transportation network
under consideration has a spectral radius $\rho\left(A_{\mathcal{G}}\right)=9.46$.
In our simulations, we consider the following bounds for the feasible
infection and recovery rates: $\underline{\delta}_{i}=0.1$, $\bar{\delta}_{i}=0.5$
and $\underline{\beta}_{i}=2.1\times10^{-2}$, $\bar{\beta}_{i}=4.2\times10^{-3}$,
for $i=1,\ldots,56$. Notice that, in the absence of protection resources,
the matrix $BA_{\mathcal{G}}-D$ in (\ref{eq:Spectral Control}) has
its largest eigenvalue at $\lambda_{1}\left(\bar{\beta}_{i}A_{\mathcal{G}}-\underline{\delta}_{i}I\right)=0.1>0$;
thus, the disease-free equilibrium is unstable. We now find the cost-optimal
allocation of resources to stabilize the disease-free equilibrium.

\begin{figure}
\centering\includegraphics[width=0.95\columnwidth]{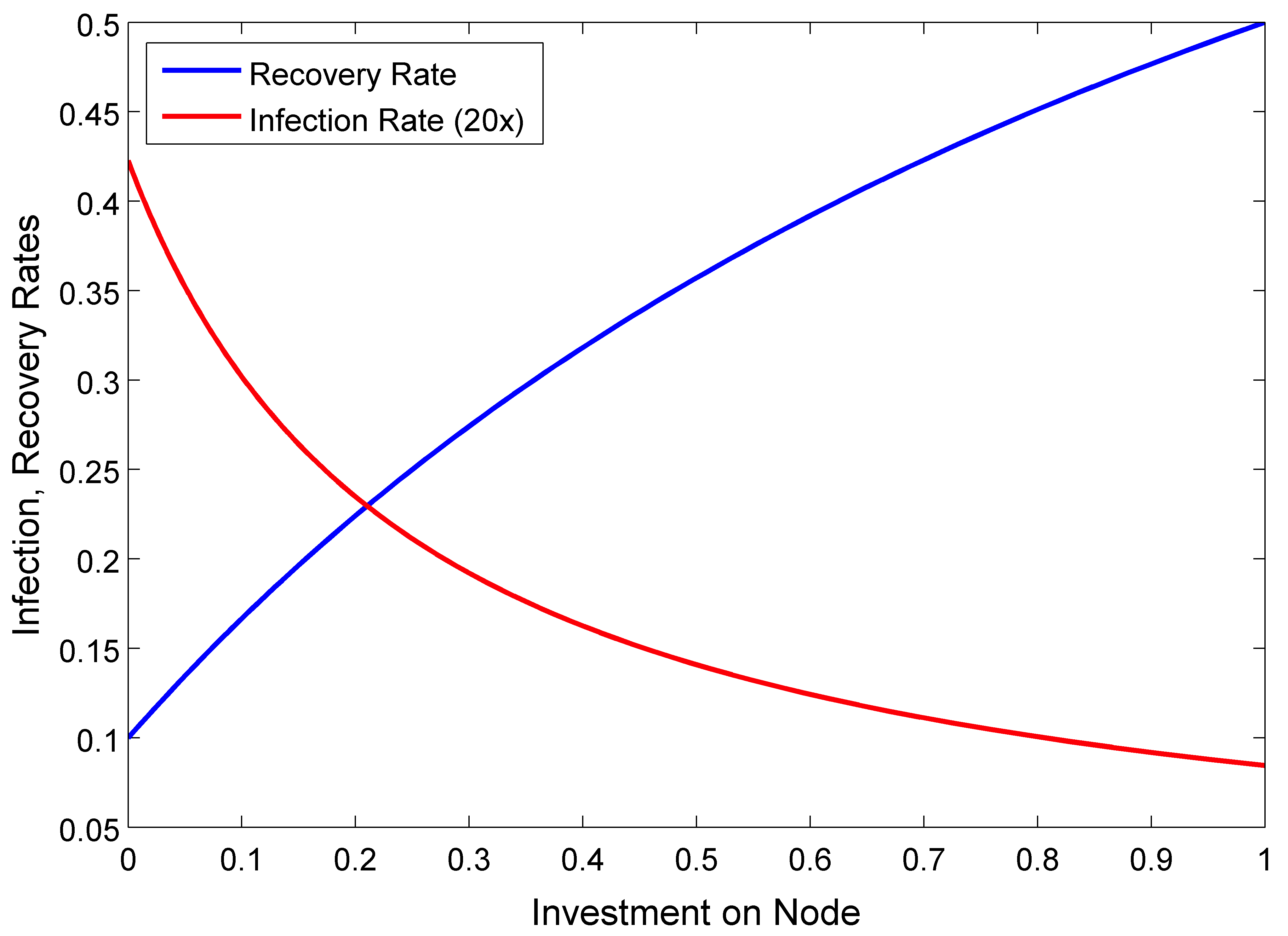}
\caption{Infection rate (in red, and multiplied by 20, to improve visualization)
and recovery rate (in blue) achieved at node $v_{i}$ after an investment
on protection (in abscissas) is made on that node.}
\label{fig. cost functions}

\end{figure}

In our simulations, we use cost functions inspired by the shape of the prototypical cost functions representing probability of failure vs. investment in systems reliability (see, for example, \cite{NGS07}). These cost functions are usually quasiconvex and present diminishing returns. In our context, the infection rate plays a role similar to the probability of failure in systems reliability. Consequently, we have chosen cost functions presenting these two features in our illustrations. In particular, we consider the following cost functions:
\begin{equation}
f_{i}\left(\beta_{i}\right)=\frac{\beta_{i}^{-1}-\bar{\beta}_{i}^{-1}}{\underline{\beta}_{i}^{-1}-\bar{\beta}_{i}^{-1}},\: g_{i}\left(\delta_{i}\right)=\frac{\left(1-\delta_{i}\right)^{-1}-\left(1-\underline{\delta}_{i}\right)^{-1}}{\left(1-\overline{\delta}_{i}\right)^{-1}-\left(1-\underline{\delta}_{i}\right)^{-1}}.\label{eq:Quasiconvex Limit-1}
\end{equation}
Notice that we have normalized these cost functions to have values
in the interval $\left[0,1\right]$ when $\underline{\beta}_{i}\leq\beta_{i}\leq\bar{\beta}_{i}$
and $\underline{\delta}_{i}\leq\delta_{i}\leq\bar{\delta}_{i}$. In
Fig. \ref{fig. cost functions}, we plot these cost functions, where
the abscissa is the amount invested in either vaccines or antidotes
on a particular node and the ordinates are the infection (red line)
and recovery (blue line) rates achieved by the investment. As we increase
the amount invested on vaccines from 0 to 1, the infection rate of
that node is reduced from $\bar{\beta}_{i}$ to $\underline{\beta}_{i}$
(red line). Similarly, as we increase the amount invested on antidotes
at a node $v_{i}$, the recovery rate grows from $\underline{\delta}_{i}$
to $\bar{\delta}_{i}$ (blue line). Notice that both cost functions
present diminishing marginal benefit on investment.

\textbf{}%

\begin{figure*}
\centering\includegraphics[width=1\textwidth]{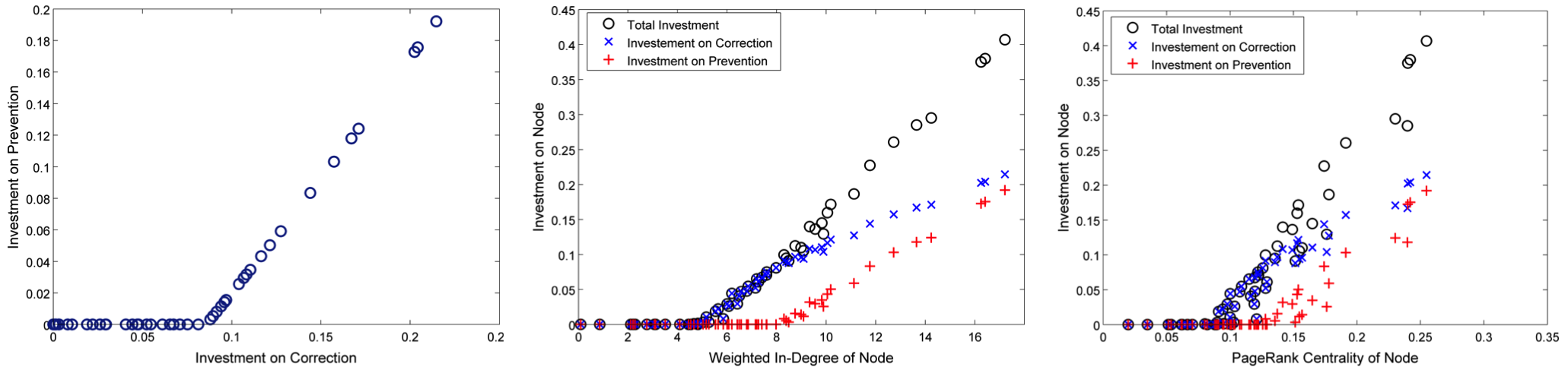}
\caption{Results from the rate-constrained allocation problem. From left to
right, we have (\emph{a}) a scatter plot with the investment on correction
versus prevention per node, (\emph{b}) a scatter plot with the investment
on protection per node and the in-degrees, and (\emph{c}) a scatter
plot with the investment on protection per node versus PageRank centralities.}

\label{fig_rate-constrained} 
\end{figure*}

Using the air transformation network, the parameters, and the cost
functions specified above, we solve both the rate-constrained and
the budget-constrained allocation problem using the geometric programs
in Theorems \ref{thm:GP for budget constrained} and \ref{thm:GP for rate constrained}.
The solution of the rate-constrained problem with $\overline{\varepsilon}=10^{-3}$
is summarized in Fig. \ref{fig_rate-constrained}. In the left subplot,
we present a scatter plot with 56 circles (one circle per airport),
where the abscissa of each circle is equal to $g_{i}\left(\delta_{i}^{*}\right)$
and the ordinate is $f_{i}\left(\beta_{i}^{*}\right)$, namely, the
investments on correction and prevention on the airport at node $v_{i}$,
respectively. We observe an interesting pattern in the allocation
of preventive and corrective resources in the network. In particular,
we have that in the optimal allocation some airports receive no resources
at all (the circles associated to those airports are at the origin
of the scatter plot); some airports receive only corrective resources
(indicated by circles located on top of the $x$-axis), and some airports
receive a mixture of preventive and corrective resources. In the center
and right subplots in Fig. \ref{fig_rate-constrained}, we compare
the distribution of resources with the in-degree and the PageRank%
\footnote{The PageRank vector $\mathbf{r}$, before normalization, can be computed
as $\mathbf{r}=(I-\alpha A_{\mathcal{G}}\mbox{diag}(1/\deg_{out}\left(v_{i}\right)))^{-1}\mathbf{1}$,
where $\mathbf{1}$ is the vector of all ones and $\alpha$ is typically
chosen to be $0.85$.%
} centralities of the nodes in the network \cite{New10}. In the center
subplot, we have a scatter plots where the ordinates represent investments
on prevention (red +'s), correction (blue x's), and total investment
(the sum of prevention and correction investments, in black circles)
for each airport, while the abscissas are the (weighted) in-degrees%
\footnote{It is worth remarking that the in-degree in the abscissa of Fig. \ref{fig_rate-constrained}
accounts from the incoming traffic into airport $v_{i}$ coming only
from those airports in the selective group of airports with an incoming
traffic over 10 MPPY. Therefore, the in-degree does not represent
the total incoming traffic into the airport.%
} of the airports under consideration. We again observe a nontrivial
pattern in the allocation of investments for protections. In particular,
for airports with incoming traffic less than $5$ MPPY, no resources
are invested at all, while for airports with incoming traffic between
$5$ MPPY and $8$ MPPY, only corrective resources are needed. Airports
with incoming traffic over $8$ MPPY receive both preventive and corrective
resources. In the right subplot in Fig. \ref{fig_rate-constrained},
we include a scatter plot of the amount invested on prevention and
correction for each airport versus its PageRank centrality in the
transportation network. We observe that, although there is a strong
correlation between centrality and investments, there is no trivial
law to achieve the optimal resource allocation based on centrality
measurements solely. For example, we observe in Fig. \ref{fig_rate-constrained}
how some airports receive a higher investment on protection than other
airports with higher centrality. Hence, it is not always best for the most central nodes to receive the most resources. In \cite{ZP14}, the authors expand on this phenomenon and propose digraphs in which resources allocated in the most central nodes have no effect on the exponential decay rate; these resources are effectively wasted.

\begin{figure*}
\centering\includegraphics[width=1\textwidth]{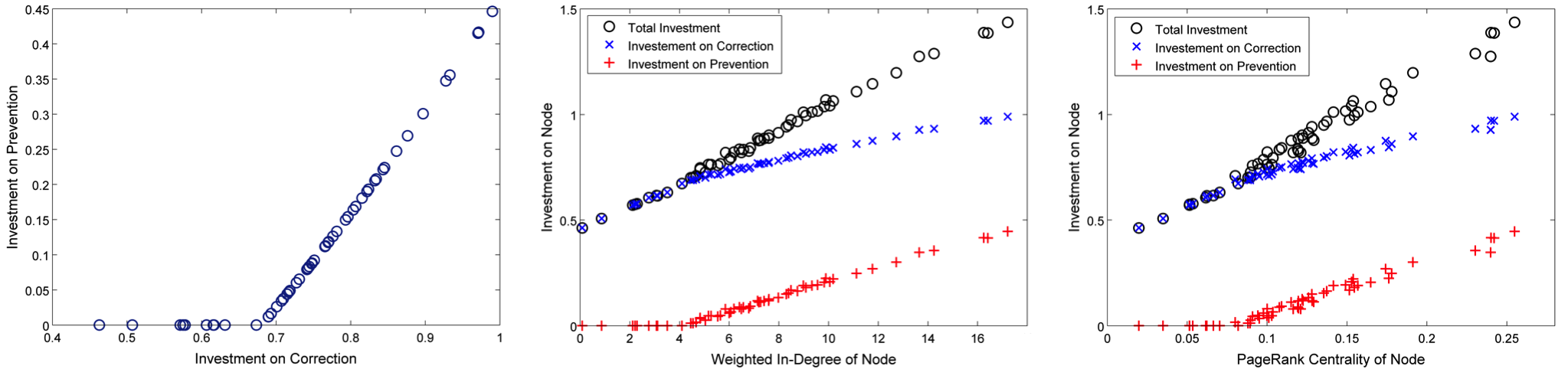}
\caption{Results from the budget-constrained allocation problem. From left
to right, we have (\emph{a}) a scatter plot with the investment on
correction versus prevention per node, (\emph{b}) a scatter plot with
the investment on protection per node and the in-degrees, and (\emph{c})
a scatter plot with the investment on protection per node versus PageRank
centralities.}

\label{fig_budget-constrained} 
\end{figure*}

Using Theorem \ref{thm:GP for budget constrained}, we also solve
the budget-constrained allocation problem. We have chosen a budget
that is a 50\% extra over the optimal budget computed from Problem
\ref{Problem:Rate Constrained Allocation} with $\overline{\varepsilon}=10^{-3}$.
With this extra budget, we achieve an exponential decay rate of $\varepsilon^{*}=0.342$.
The corresponding allocation of resources is summarized in Fig. \ref{fig_budget-constrained}.
The subplots in this figure are similar to those in Fig. \ref{fig_rate-constrained},
and we only remark the main differences in here. Notice that, given
the extra budget, there are no airports with no investment on protection
resources (as indicated by the absence of circles at the origin of
the left subplot). Also, the center subplot indicates that all the
airports receive a certain amount of corrective resources, although
not all of them receive preventive resources (such as those with a
(weighted) in-degree less than 4 MPPY). The scatter plot at the right
illustrates the relationship between investments and PageRank centrality.

{}

\section{\label{sec:Conclusions}Conclusions}

We have studied the problem of allocating protection resources in
weighted, directed networks to contain spreading processes, such as
the propagation of viruses in computer networks, cascading failures
in complex technological networks, or the spreading of an epidemic
in a human population. We have considered two types of protection
resources: (\emph{i}) \emph{Preventive} resources able to `immunize'
nodes against the spreading (e.g. vaccines), and (\emph{ii}) \emph{corrective}
resources able to neutralize the spreading after it has reached a
node (e.g. antidotes). We assume that protection resources have an
associated cost and have then studied two optimization problems: (\emph{a})
The \emph{budget-constrained allocation problem}, in which we find
the optimal allocation of resources to contain the spreading given
a fixed budget, and (\emph{b}) the \emph{rate-constrained allocation
problem}, in which we find the cost-optimal allocation of protection
resources to achieve a desired decay rate in the number of `infected'
nodes.

We have solved these
optimal resource allocation problem in \emph{weighted and directed}
networks of \emph{nonidentical} agents in polynomial time using Geometric
Programming (GP). Furthermore, the framework herein proposed allows
\emph{simultaneous }optimization over both preventive and corrective
resources, even in the case of cost functions being node-dependent.

We have illustrated our approach by designing an optimal protection
strategy for a real air transportation network. We have limited our
study to the network of the world's busiest airports by passenger
traffic. For this transportation network, we have computed the optimal
distribution of protecting resources to contain the spread of a hypothetical
world-wide pandemic. Our simulations show that the optimal distribution
of protecting resources follows nontrivial patterns that cannot, in
general, be described using simple heuristics based on traditional
network centrality measures.

\appendix{}

\textbf{Proof of Lemma \ref{lem:Monotonicity of lambda1}.} We define
the auxiliary matrix $M\triangleq\mbox{diag}\left(\beta_{i}\right)A-\mbox{diag}\left(\delta_{i}\right)+\Delta I$,
where $\Delta\triangleq\max\left\{ \delta_{i}\right\} $. Thus, $\lambda_{1}\left(M\right)=\lambda_{1}\left(\mbox{diag}\left(\beta_{i}\right)A-\mbox{diag}\left(\delta_{i}\right)\right)+\Delta$.
Notice that both $M$ and $M^{T}$ are nonnegative and irreducible
if $\mathcal{G}$ is strongly connected. Hence, from Lemma \ref{lem:Perron-Frobenius},
there are two positive vectors $\mathbf{v}$ and $\mathbf{w}$ such
that
\begin{align*}
M\mathbf{v} & =\rho\mathbf{v},\\
\mathbf{w}^{T}M & =\rho\mathbf{w}^{T},
\end{align*}
where $\rho=\rho\left(M\right)=\lambda_{1}\left(M\right)$, and $\mathbf{v}$,
$\mathbf{w}$ are the right and left dominant eigenvectors of $M$.
From eigenvalue perturbation theory, we have that the increment in
the spectral radius of $M$ induced by a matrix increment $\Delta M$
is \cite{meyer2000matrix}
\begin{equation}
\rho\left(M+\Delta M\right)-\rho\left(M\right)=\mathbf{w}^{T}\Delta M\mathbf{v}+o\left(\left\Vert \Delta M\right\Vert \right).\label{eq:Eigenvalue Perturbation}
\end{equation}

To study the effect of a positive increment in $\beta_{k}$ in the
spectral radius, we define $\Delta B=\Delta\beta_{k}\mathbf{e}_{k}\mathbf{e}_{k}^{T}$,
for $\Delta\beta_{k}>0$, and apply \ref{eq:Eigenvalue Perturbation}
with $\Delta M=\Delta B\, A$. Hence,
\begin{align*}
\rho\left(M+\Delta M\right)-\rho\left(M\right) & =\Delta\beta_{k}\mathbf{w}^{T}\mathbf{e}_{k}\mathbf{e}_{k}^{T}A\mathbf{v}+o\left(\left\Vert \Delta\beta_{k}\right\Vert \right)\\
 & =\Delta\beta_{k}w_{k}\mathbf{a}_{k}^{T}\mathbf{v}+o\left(\left\Vert \Delta\beta_{k}\right\Vert \right)>0,
\end{align*}
where $\mathbf{a}_{k}^{T}=\mathbf{e}_{k}^{T}A$ and the last inequality
if a consequence of $\Delta\beta_{k}$, $w_{k}$, and $\mathbf{a}_{k}^{T}\mathbf{v}$
being all positive. Hence, a positive increment in $\beta_{k}$ induces
a positive increment in the spectral radius.

\begin{flushright}
Similarly, to study the effect of a positive increment in $\delta_{k}$
in the spectral radius, we define $\Delta D=\Delta\delta_{k}\mathbf{e}_{k}\mathbf{e}_{k}^{T}$,
for $\Delta\delta_{k}>0$. Applying \ref{eq:Eigenvalue Perturbation}
with $\Delta M=-\Delta D$, we obtain
\begin{align*}
\rho\left(M+\Delta D\right)-\rho\left(M\right) & =-\Delta\delta_{k}\mathbf{w}^{T}\mathbf{e}_{k}\mathbf{e}_{k}^{T}\mathbf{v}+o\left(\left\Vert \Delta\delta_{k}\right\Vert \right)\\
 & =-\Delta\delta_{k}w_{k}v_{k}+o\left(\left\Vert \Delta\delta_{k}\right\Vert \right)<0.
\end{align*}
$\blacksquare$
\par\end{flushright}

\textbf{Proof of Lemma \ref{lem:Pattern of Zeros}.} The proof of
(\emph{a}) is trivial and valid for any square matrix $M$. To prove
(\emph{b}), we consider the eigenvalue equations for $M$ and $RM$,
i.e., $M\mathbf{u}=\lambda_{1}\left(M\right)\mathbf{u}$ and $RM\mathbf{w}=\lambda\mathbf{w}$,
where $\mathbf{u}=\mathbf{v}_{1}\left(M\right)=\left[u_{i}\right]$
and $\mathbf{w}=\mathbf{v}_{1}\left(RM\right)=\left[w_{i}\right]$.
We expand the eigenvalue equations component-wise as,
\begin{align}
\sum_{j=1}^{n}m_{ij}u_{j} & =\lambda u_{i},\label{eq:Eig M}\\
\sum_{j=1}^{n}r_{i}m_{ij}w_{j} & =\lambda w_{i},\label{eq:Eig RM}
\end{align}
for all $i=1,\ldots,n$. We now prove statement (\emph{b}) by proving
that $v_{i}=0$ if and only if $w_{i}=0$.

If $u_{i}=0$, then (\ref{eq:Eig M}) gives $\sum_{j}m_{ij}v_{j}=0$.
Since $m_{ij},v_{i}\geq0$, the summation $\sum_{j}m_{ij}v_{j}=0$
if and only if the following two statements hold: (\emph{a1}) $m_{ij}>0\implies v_{j}=0$
and (\emph{a2}) $v_{j}>0\implies m_{ij}=0$. Since $t_{i}>0$, these
two statements are equivalent to: (\emph{b1}) $t_{i}m_{ij}>0\implies v_{j}=0$
and (\emph{b2}) $v_{j}>0\implies t_{i}m_{ij}=0$. Statements (\emph{b1})
and (\emph{b2}) are true if and only if $\sum_{j}\left(t_{i}m_{ij}\right)w_{j}=0=w_{i}$,
where the last equality comes from (\ref{eq:Eig RM}). Hence, we have
that $v_{i}=0\iff w_{i}=0$; hence, $\mathcal{Z}\left(\mathbf{u}\right)=\mathcal{Z}\left(\mathbf{w}\right)$.

\begin{flushright}
$\blacksquare$
\par\end{flushright}

\textbf{Proof of Proposition \ref{prop:Zeros in BA-D}.} Our proof
is based on the transformations defined in Lemma \ref{lem:Pattern of Zeros}.
Starting from a matrix $BA-D$, we then apply the following chain
of transformations:
\begin{description}
\item [{\textmd{(}\textmd{\emph{i}}\textmd{)}}] $T_{\alpha}\left(BA-D\right)=BA+\Delta$,
for $\alpha=\max\left\{ d_{i}\right\} $. Hence, $\Delta=\max\left\{ d_{i}\right\} I-D$
and $BA+\Delta\geq0$.
\item [{\textmd{(}\textmd{\emph{ii}}\textmd{)}}] $T_{R}\left(BA+\Delta\right)=\Delta^{-1}BA+I$,
for $R=\Delta^{-1}$.
\item [{\textmd{(}\textmd{\emph{iii}}\textmd{)}}] $T_{\alpha}\left(\Delta^{-1}BA+I\right)=\Delta^{-1}BA$,
for $\alpha=-1$.
\item [{\textmd{(}\textmd{\emph{iv}}\textmd{)}}] $T_{R}\left(\Delta^{-1}BA\right)=A$,
for $R=B^{-1}\Delta$.
\end{description}
From Lemma \ref{lem:Pattern of Zeros}, these transformations preserve
the location of the zeros in the dominant eigenvector. Thus, the input
to the first transformation, $BA-D$, and the output to the last transformation,
$A$, satisfy $\mathcal{Z}\left(\mathbf{v}_{1}\left(BA-D\right)\right)=\mathcal{Z}\left(\mathbf{v}_{1}\left(A\right)\right)$.

\begin{flushright}
$\blacksquare$
\par\end{flushright}

\bibliographystyle{ieeetr}
\bibliography{ViralSpread}

\begin{thebibliography}{10}

\bibitem{Bai75}
N.~Bailey, {\em The mathematical theory of infectious diseases and its
  applications}.
\newblock Charles Griffin \& Company Ltd., 1975.

\bibitem{GGT03}
M.~Garetto, W.~Gong, and D.~Towsley, ``Modeling malware spreading dynamics,''
  in {\em IEEE INFOCOM 2003}, vol.~3, pp.~1869--1879, 2003.

\bibitem{roy2012security}
S.~Roy, M.~Xue, and S.~Das, ``Security and discoverability of spread dynamics
  in cyber-physical networks,'' {\em IEEE Transactions on Parallel and
  Distributed Systems}, vol.~23, no.~9, pp.~1694--1707, 2012.

\bibitem{WCWF03}
Y.~Wang, D.~Chakrabarti, C.~Wang, and C.~Faloutsos, ``Epidemic spreading in
  real networks: An eigenvalue viewpoint,'' in {\em Proc. 22nd Int. Symp.
  Reliable Distributed Systems}, pp.~25--34, 2003.

\bibitem{GMT05}
A.~Ganesh, L.~Massoulie, and D.~Towsley, ``The effect of network topology on
  the spread of epidemics,'' in {\em IEEE INFOCOM 2005}, vol.~2,
  pp.~1455--1466, 2005.

\bibitem{MOK09}
P.~Van~Mieghem, J.~Omic, and R.~Kooij, ``Virus spread in networks,'' {\em
  IEEE/ACM Transactions on Networking}, vol.~17, no.~1, pp.~1--14, 2009.

\bibitem{cohen2003efficient}
R.~Cohen, S.~Havlin, and D.~Ben-Avraham, ``Efficient immunization strategies
  for computer networks and populations,'' {\em Physical Review Letters},
  vol.~91, no.~24, p.~247901, 2003.

\bibitem{BCGS10}
C.~Borgs, J.~Chayes, A.~Ganesh, and A.~Saberi, ``How to distribute antidote to
  control epidemics,'' {\em Random Structures \& Algorithms}, vol.~37, no.~2,
  pp.~204--222, 2010.

\bibitem{chung2009distributing}
F.~Chung, P.~Horn, and A.~Tsiatas, ``Distributing antidote using pagerank
  vectors,'' {\em Internet Mathematics}, vol.~6, no.~2, pp.~237--254, 2009.

\bibitem{WRS08}
Y.~Wan, S.~Roy, and A.~Saberi, ``Designing spatially heterogeneous strategies
  for control of virus spread,'' {\em IET Systems Biology}, vol.~2, no.~4,
  pp.~184--201, 2008.

\bibitem{GOM11}
E.~Gourdin, J.~Omic, and P.~Van~Mieghem, ``Optimization of network protection
  against virus spread,'' in {\em 8th International Workshop on the Design of
  Reliable Communication Networks}, pp.~86--93, 2011.

\bibitem{PDS13}
V.~Preciado, F.~Darabi~Sahneh, and C.~Scoglio, ``A convex framework for optimal
  investment on disease awareness in social networks,'' in {\em IEEE Global
  Conference on Signal and Information Processing}, pp.~851--854, 2013.

\bibitem{PZ13}
V.~Preciado and M.~Zargham, ``Traffic optimization to control epidemic
  outbreaks in metapopulation models,'' in {\em IEEE Global Conference on
  Signal and Information Processing}, pp.~847--850, 2013.

\bibitem{PZEJP13}
V.~Preciado, M.~Zargham, C.~Enyioha, A.~Jadbabaie, and G.~Pappas, ``Optimal
  vaccine allocation to control epidemic outbreaks in arbitrary networks,'' in
  {\em IEEE Conference on Decision and Control}, 2013.

\bibitem{weiss1971asymptotic}
G.~Weiss and M.~Dishon, ``On the asymptotic behavior of the stochastic and
  deterministic models of an epidemic,'' {\em Mathematical Biosciences},
  vol.~11, no.~3, pp.~261--265, 1971.

\bibitem{VO13}
P.~Van~Mieghem and J.~Omic, ``In-homogeneous virus spread in networks,'' {\em
  arXiv preprint arXiv:1306.2588}, 2013.

\bibitem{barrat2008dynamical}
A.~Barrat, M.~Barthelemy, and A.~Vespignani, {\em Dynamical processes on
  complex networks}.
\newblock Cambridge University Press Cambridge, 2008.

\bibitem{GMWPM12}
J.~Gleeson, S.~Melnik, J.~Ward, M.~Porter, and P.~Mucha, ``Accuracy of
  mean-field theory for dynamics on real-world networks,'' {\em Physical Review
  E}, vol.~85, p.~026106, 2012.

\bibitem{van2006performance}
P.~Van~Mieghem, {\em Performance analysis of communications networks and
  systems}.
\newblock Cambridge University Press, 2006.

\bibitem{BV04}
S.~Boyd and L.~Vandenberghe, {\em Convex optimization}.
\newblock Cambridge university press, 2004.

\bibitem{BKVH07}
S.~Boyd, S.-J. Kim, L.~Vandenberghe, and A.~Hassibi, ``A tutorial on geometric
  programming,'' {\em Optimization and engineering}, vol.~8, no.~1,
  pp.~67--127, 2007.

\bibitem{meyer2000matrix}
C.~Meyer, {\em Matrix analysis and applied linear algebra}.
\newblock SIAM, 2000.

\bibitem{New10}
M.~Newman, {\em Networks: An introduction}.
\newblock Cambridge University Press, 2010.

\bibitem{schneider2011suppressing}
C.~Schneider, T.~Mihaljev, S.~Havlin, and H.~Herrmann, ``Suppressing epidemics
  with a limited amount of immunization units,'' {\em Physical Review E},
  vol.~84, no.~6, p.~061911, 2011.

\bibitem{NGS07}
E.~Nikolaidis, D.~Ghiocel, and S.~Singhal, {\em Engineering Design Reliability
  Applications: For the Aerospace, Automotive and Ship Industries}.
\newblock CRC Press, 2007.

\bibitem{ZP14}
M.~Zargham and V.~Preciado, ``Worst-case scenarios for greedy, centrality-based
  network protection strategies,'' in {\em ArXiv:1401.5753v1}, 2014.

\end{thebibliography}

\end{document}